\documentclass[12pt]{article}
\usepackage{amssymb, amsmath}
\begin{document}
\voffset=0.0truein \hoffset=-0.5truein \setlength{\textwidth}{6.0in}
\setlength{\textheight}{8.8in} \setlength{\topmargin}{-0.2in}
\renewcommand{\theequation}{\arabic{section}.\arabic{equation}}
\newtheorem{theorem}{Theorem}[section]
\newtheorem{lemma}{Lemma}[section]
\newtheorem{pro}{Proposition}[section]
\newtheorem{cor}{Corollary}[section]
\newcommand{\n}{\nonumber}
\newcommand{\e}{\varepsilon}
\renewcommand{\a}{\alpha}
\renewcommand{\l}{\lambda}
\newcommand{\lo}{\lambda_1}
\newcommand{\lt}{\lambda_2}
\newcommand{\ltt}{\lambda_3}
\newcommand{\s}{\sigma}
\renewcommand{\o}{\omega}
\renewcommand{\O}{\Bbb T^3}
\newcommand{\Om}{\Omega}
\newcommand{\intt}{\int_{\Bbb T^3}}
\newcommand{\intr}{\int_{\Bbb T^3}}
\newcommand{\bb}{\begin{equation}}
\newcommand{\ee}{\end{equation}}
\newcommand{\bq}{\begin{eqnarray}}
\newcommand{\eq}{\end{eqnarray}}
\newcommand{\bqn}{\begin{eqnarray*}}
\newcommand{\eqn}{\end{eqnarray*}}
\title{On the spectral dynamics of the deformation tensor and new a
priori estimates for  \\
the 3D Euler equations}
\author{Dongho Chae\thanks{Permanent Address: Department of Mathematics,
              Sungkyunkwan University, Suwon 440-746, Korea,
              e-mail : chae@skku.edu}\\
Center for Scientific Computation And\\
 Mathematical Modeling\\
 Paint Branch Drive \\
 University of Maryland\\
  College Park, MD 20742-3289, USA\\
  e-mail: {\it dchae@cscamm.umd.edu}}
 \date{}
\maketitle
\begin{abstract}
In this paper we study the dynamics of eigenvalues of the
deformation tensor for  solutions of the 3D incompressible Euler
equations.   Using the evolution equation of the $L^2$ norm of
spectra, we deduce new a priori estimates of the $L^2$ norm of
vorticity. As an immediate corollary of the estimate we obtain a new
sufficient condition of $L^2$ norm control of vorticity. We also
obtain decay in time estimates of the ratios of the eigenvalues. In
the remarks we discuss  what these estimates suggest in the study of
searching initial data leading to a possible finite time
singularities. We find that the dynamical behaviors of $L^2$ norm of
vorticity are controlled completely by the second largest
eigenvalue of the deformation tensor.
\end{abstract}

\section{Introduction}
\setcounter{equation}{0}

 We are  concerned with the following Euler
equations for the homogeneous incompressible fluid flows in $\Omega
\subset \Bbb R^3$.
 \bq
 \frac{\partial v}{\partial t} +(v\cdot \nabla )v &=-\nabla p ,
  \\
 \textrm{div }\, v =0 , &\\
  v(x,0)=v_0 (x), &
  \eq
where $v=(v_1, v_2, v_3 )$, $v_j =v_j (x, t)$, $j=1, 2,3, $ is the
velocity of the flow, $p=p(x,t)$ is the scalar pressure, and $v_0 $
is the given initial velocity, satisfying div $v_0 =0$. For
simplicity of presentation we assume $\Omega = \Bbb T^3$, the 3D
periodic box. Most of the results in this paper, however, are valid
also in the whole of $\Bbb R^3$, or bounded domain with smooth
boundary, at least after obvious modifications. Given $m\in \Bbb
N\cup\{0\}$, let  $H^m (\O)$ be the standard Sobolev space on $\Bbb
T^3$,
$$H^m (\Bbb T^3)=\{ f\in L^2 (\Bbb T^3 )\,|\,
\|f\|_{H^m}^2 =\sum_{|\alpha| \leq m}\intr |D^\alpha f(x)|^2 dx
<\infty\},$$
 where $\alpha=(\alpha_1,\alpha_2,\alpha_3 )$ with
 $|\alpha|=\alpha_1+\alpha_2+\alpha_3$ is the usual
 multi-index notation.
 We introduce the space of solenoidal vector fields,
$$\mathbb{H}^m _\sigma =\{ u \in [H^m (\Bbb T^3 )]^3\, |\,
\mathrm{div}\, u=0\}.$$
 Then, for  $ v_0 \in  \mathbb{H}^m _\sigma$ with  $m>5/2$,
  the local in
time unique existence of solution  to (1.1)-(1.3), which belongs to
$C([0,T];\mathbb{H}^m _\sigma )$ for some $T=T(\|v_0\|_{H^m})$, was
established in \cite{kat, tem}. This was later extended to the local
existence in various other function spaces by many
authors(\cite{kat1,che1,che2, vis1,vis2, cha1, cha2, cha3,cha4}).
The question of finite time blow-up/global regularity of such
locally constructed solution is one of the most outstanding open
problems in the mathematical fluid mechanics. For physical meaning
and other significance of this problem as well as many instructive
examples of solutions we refer \cite{con1,maj}. For a mathematical
or numerical test of the actual finite time blow-up of  a given
solution, it is important to have a good blow-up criterion.  In this
direction there is a celebrated result by
Beale-Kato-Majda(\cite{bea}), now called the BKM criterion. This
criterion is later refined in \cite{koz,cha4, che2}, using refined
versions of logarithmic Sobolev inequalities. As for another
approach to the blow-up criterion, there is a pioneering work on the
geometric type of blow-up criterion due to
Constantin-Fefferman-Majda(\cite{con2})(see also \cite{con1}), which
was initiated in \cite{con4}, and  the idea of which was used in the
recent works in \cite{cha5, den}. For different type of geometric
approach to the blow-up problem(considering vortex tube initial
data) we refer \cite{cor}. We also mention a recent interesting
result in \cite{bab} for rotating flows, where they discovered that
rotation has some sense of regularization effect. In this paper we
study the regularity/blow-up problem, using the spectral dynamics of
the deformation tensor for the solution of the Euler equations.
Previous spectral approaches to the singularity problems in the
nonlinear partial differential equations are studied  in \cite{liu},
however their study is not for the real Euler equations, but for its
model equations, avoiding the difficulty of  the nonlocality caused
by the Riesz transform appearing in the equations when the pressure
is eliminated. Moreover, their spectrum is for the matrices of the
velocity gradient, not for the deformation tensor, which is the
symmetric part of the velocity gradient. In the next section we
derive an evolution equation, in the $L^2$ sense, of the eigenvalues
of the deformation tensor. From this equation we derive new a priori
estimates for the $L^2$ norm of vorticity.  The inequality itself
 already tells us interesting dynamical mechanism of compression and stretching
of infinitesimal fluid volume elements leading to possible blow-up.
The inequality immediately leads to very simple and elegant
sufficient condition of $L^2$ norm control of vorticity of the 3D
Euler equations. In the section 3 we consider special classes of
initial data. For such initial data we can have better estimates
exponential growth/decay of the $L^2$ norm of vorticity. We also
deduce decay estimates in time
of a ratio of eigenvalues of the deformation tensor.\\

\section{Dynamics of eigenvalues of the deformation tensor}
\setcounter{equation}{0}

We use the following notations for matrix components.
$$
V_{ij}=\frac{\partial v_j}{\partial x_i},\quad
S_{ij}=\frac{V_{ij}+V_{ji}}{2},\quad A_{ij}=\frac{V_{ij}-V_{ji}}{2},
$$
where $i,j=1,2,3$. Then, obviously we have $V_{ij}=S_{ij}+A_{ij}$.
For a given 3D velocity field $v(x,t)$, describing fluid motions,
$S_{ij}$ is called the deformation tensor, while $A_{ij}$ is related
to the vorticity field $\o=$ curl $v$ by
$$ A_{ij}=\frac12 \sum_{k=1}^3\e_{ijk} \o_k , $$
 where $\e_{ijk}$, the skewsymmetric tensor with the normalization
 $\e_{123}=1$. For incompressible fluid we have
 $
  Tr (S)=\sum_{i=1}^3 S_{ii} =\mathrm{div}\, v=0.
  $
 We now state the theorem on the evolutions of the eigenvalues of
 the deformation tensor associated with the solution of the Euler
 system (1.1)-(1.3).
\begin{theorem}
 Let $\lo (x,t),\lt (x,t) , \ltt (x,t)$ be the eigenvalues of the deformation tensor
 $S=(S_{ij})_{i,j=1}^3$ associated to the classical solution
 $v(x,t)$ of (1.1)-(1.3). Then,  the following equation holds.
 \bb\label{2.1}
 \frac{d}{dt} \intr (\lo^2 +\lt^2 +\ltt^2 )dx=-4\intr \lo\lt\ltt dx.
 \ee
 \end{theorem}
{\bf Proof.} We take $L^2$ inner product (1.1) with $\Delta u$, and
integrate by
 part to derive
 \bq\label{2.2}
 \lefteqn{\frac12 \frac{d}{dt} \|\nabla v\|_{L^2}^2 =\intr (v\cdot \nabla )v
 \cdot \Delta v dx
 =-\sum_{i,j,k=1}^3\intr \frac{\partial v_j}{\partial x_k}
 \frac{\partial v_k}{\partial x_i}\frac{\partial v_j}{\partial x_i}
 dx }\hspace{.1in}\n \\
&=&-\sum_{i,j,k=1}^3\intr S_{kj} V_{ik} V_{ij} dx
=-\sum_{i,j,k=1}^3\intr S_{kj} (S_{ik}+ A_{ik}) (S_{ij}+ A_{ij})
dx\n \\
 &=&-\sum_{i,j,k=1}^3\intr ( S_{kj} A_{ik} A_{ij}
+S_{kj}S_{ik}S_{ij}
)dx\n \\
&=&- \frac14\sum_{i,j,k=1}^3\intr S_{kj}\left[
\sum_{m=1}^3\delta_{kj}\o_m\o_m -\o_j\o_k \right]
dx-\sum_{i,j,k=1}^3\intr S_{kj}S_{ik}S_{ij}
dx\n \\
&=&\frac14 \sum_{j,k=1}^3 \intr S_{jk} \o_j \o_k dx
-\sum_{i,j,k=1}^3\intr S_{kj}S_{ik}S_{ij} dx.
 \eq
 Next, we consider the vorticity equation for the 3D Euler equations,
 \bb\label{2.3}
 \frac{\partial \o}{\partial t} +(v \cdot \nabla )\o =(\o \cdot
 \nabla )v ,
 \ee
 which is obtained from (1.1) by taking curl$(\cdot)$ operation.
 Taking $L^2$ inner product (\ref{2.3}) with $\o$, we obtain, after
 integration by part,
 \bb\label{2.4}
  \frac12 \frac{d}{dt} \|\o \|_{L^2}^2 =\intr (\o \cdot \nabla )v
  \cdot \o dx =\sum_{j,k=1}^3\intr S_{jk} \o_j \o_k dx .
  \ee
  Since we have the equality,
  \bb\label{2.5}
  \intr |\nabla v |^2 dx =\intr |\o |^2 dx,
  \ee
  from (\ref{2.2}) and (\ref{2.4}) we obtain
  $$
\frac14 \sum_{j,k=1}^3 \intr S_{jk} \o_j \o_k dx
-\sum_{i,j,k=1}^3\intr S_{kj}S_{ik}S_{ij} dx=\sum_{j,k=1}^3\intr
S_{jk} \o_j \o_k dx, $$
 Hence,
\bb\label{2.6}
 \sum_{j,k=1}^3\intr
S_{jk} \o_j \o_k dx=-\frac43 \sum_{i,j,k=1}^3\intr
S_{kj}S_{ik}S_{ij} dx. \ee
 We also have the following pointwise equality,
 \bq\label{2.7}
 |\nabla v|^2&=&\sum_{j,k=1}^3 V_{jk}V_{jk}=\sum_{j,k=1}^3 (S_{jk}
 +A_{jk} )(S_{jk}
 +A_{jk} )\n\\
 &=&\sum_{j,k=1}^3 \left(S_{jk}S_{jk} +\frac14\sum_{m,n}^3
 \e_{jkm}\e_{jkn}\o_m\o_n \right)\n \\
 &=&\sum_{j,k=1}^3 S_{jk}S_{jk} +\frac12 |\o |^2.
 \eq
 Integrating (\ref{2.7}) over $\O$, and using  (\ref{2.5}), we
 obtain
 \bb\label{2.8}
 \intr |\o |^2 dx  =2\sum_{j,k=1}^3 \intr S_{jk}S_{jk} dx=2\intr (\lo^2 +\lt^2 +\ltt^2 )dx.
 \ee
 We also observe,
 \bb\label{2.9}
\sum_{i,j,k=1}^3 S_{kj}S_{ik}S_{ij} = \lo^3 +\lt^3 +\ltt^3=
3\lo\lt\ltt,
 \ee
 which follows from the following algebra, using $\lo+\lt+\ltt=0$,
 \bqn
 0&=&(\lo+\lt +\ltt )^3\\
 &=&\lo^3 +\lt^3 +\ltt^3 +3\lo^2 (\lt+\ltt
 )+3\lt^2 (\lo +\ltt ) +3\ltt (\lo+\lt )+ 6\lo\lt\ltt\\
 &=&\lo^3 +\lt^3 +\ltt^3 -3(\lo^3 +\lt^3 +\ltt^3)+ 6\lo\lt\ltt.
 \eqn
Substituting (2.8) and (2.6), combined with (2.9), into (2.4), we have (2.1). $\square$\\
\ \\

The following is a new a priori estimate for the $L^2$ norm of
vorticity for the 3D incompressible Euler equations.

\begin{theorem}
Let $v(t)\in C([0,T); \mathbb{H}^m _\s )$, $m>5/2$ be the local
classical solution of the 3D Euler equations with initial data $v_0
\in \mathbb{H}^m_\s $. Let $\lo (x,t)\geq \lt (x,t)\geq \ltt (x,t)$
are the eigenvalues of the deformation tensor $S_{ij}(v)=\frac12 (
\frac{\partial v_j}{\partial x_i} +\frac{\partial v_i}{\partial
x_j})$. We denote $\lt^+ (x,t)=\max\{ \lt (x,t), 0\}$, and $\lt^-
(x,t)=\min\{ \lt (x,t), 0\}$. Then, the following (a priori)
estimates hold.
 \bq\label{important}
  \lefteqn{\|\o_0\|_{L^2} \exp\left[
  \int_0 ^t \left(\frac12 \inf_{x\in \Bbb T^3}\lt^+ (x,t)-\sup_{x\in \O}
  |\lt^-(x,t)|\right)dt\right]\leq
  \|\o (t)\|_{L^2}}\hspace{1.in}\n \\
  &&\leq  \|\o_0\|_{L^2} \exp\left[
  \int_0 ^t \left( \sup_{x\in \Bbb T^3}\lt^+ (x,t)-\frac12 \inf_{x\in \O}
  |\lt^- (x,t)|\right)dt\right]\n \\
\eq
 for all $t\in (0,T)$.
\end{theorem}
{\bf Remark 2.1} The above estimate says, for example, that if we
have the following comparability conditions,
$$\sup_{x\in \Bbb T^3}\lt^+ (x,t)\simeq  \inf_{x\in \O}|\lt^- (x,t)|\simeq  g(t)$$
for some time interval $[0,T]$, then
$$ \|\o (t)\|_{L^2}
\lesssim O\left(\exp\left[C \int_0 ^t g(s)ds\right]\right) \qquad
\forall t\in [0,T]
$$
for some constant $C$.\\

 \noindent{\bf Remark 2.2} We note that
$\lt^+(x,t)>0$ implies we have stretching of infinitesimal fluid
volume in two directions and compression in the other one
direction(planar stretching) at $(x,t)$, while $|\lt^-(x,t)|>0$
implies stretching in one direction and compressions in two
directions(linear stretching). The above estimate says that the
dominance competition between planar stretching and linear
stretching is an important mechanism controlling  the  growth/decay
in time of the
$L^2$ norm of  vorticity.\\

\noindent{\bf Proof of Theorem 2.2}  Since $\lo+\lt+\ltt =0$, and
$\lo \geq \lt \geq \ltt $, we have $\lo\geq 0, \ltt \leq 0$, and
 \bb\label{basic}
 |\lt|\leq \min\{ \lo, |\ltt|\}.
 \ee
 We first
observe that from (\ref{2.8}),
 \bq\label{2.8a}
 \intr |\o |^2 dx &=&2\intr (\lo^2 +\lt^2 +\ltt^2 )dx\n\\
 &=&4\intr (\lo^2 +\lo\lt +\lt^2 )dx\quad (\ltt=-\lo-\lt )\n\\
 &=&4\intr (\lt^2 +\lt\ltt +\ltt^2 )dx\quad (\lo=-\lt-\ltt )
 \eq
We estimate the `vortex stretching term' as
 \bqn
 \lefteqn{-4\intr \lo\lt\ltt dx=
 -4\intr \lt^+\lo\ltt dx-4\intr \lt^-\lo\ltt dx}\hspace{1.in}\n \\
&=&4\intr \lt^+\lo(\lo+\lt) dx-4\intr |\lt^-|(\lt+\ltt)\ltt dx\n\\
&=&4\intr \lt^+ (\lo^2+\lo\lt) dx-2\intr |\lt^-|(2\lt\ltt +2\ltt^2) dx\n\\
&\leq&4\sup_{x\in \O}\lt^+ (x,t) \intr (\lo^2+\lo\lt+\lt^2) dx\n \\
&&\qquad -2\inf_{x\in \O}|\lt^- |(x,t)\intr (\lt^2+ \lt\ltt +\ltt^2) dx\n \\
\eqn
 \bq\label{right}
&=&4\sup_{x\in \O}\lt^+ (x,t) \intr (\lo^2+\lo\lt+\lt^2) dx\n \\
&&\qquad-2\inf_{x\in \O}|\lt^- (x,t)|\intr (\lo^2+\lo\lt+\lt^2) dx,
 \eq
 where we used (\ref{basic}) and (\ref{2.8a}). This, combined with
(\ref{2.1}) and (\ref{2.8a}), yields
 \bb
 \frac{d}{dt}\intr |\o (x,t)|^2 dx\leq \left[2\sup_{x\in \O}\lt^+ (x,t)
-\inf_{x\in \O}|\lt^- (x,t)|\right]\intr |\o (x,t)|^2 dx.
 \ee
 Applying the Gronwall lemma, we have the second inequality
 of (\ref{important}).
 In order to prove the first inequality of (\ref{important})  we
 estimate from below starting from
 equality part of (\ref{right})
\bq\label{left}
 \lefteqn{-4\intr \lo\lt\ltt dx
=4\intr \lt^+\lo(\lo+\lt) dx-4\intr |\lt^-|(\lt+\ltt)\ltt dx}\hspace{1.in}\n \\
&=&2\intr \lt^+(2\lo^2+2\lo\lt) dx-4\intr |\lt^-|(\lt\ltt +\ltt^2) dx\n\\
 &\geq&2\intr \lt^+(\lo^2+\lo\lt+\lt^2) dx-4\intr |\lt^-|(\lt^2+\lt\ltt +\ltt^2) dx\n\\
 &\geq& 2\inf_{x\in \O}\lt^+ (x,t) \intr (\lo^2+\lo\lt+\lt^2) dx\n \\
 &&\qquad-4\sup_{x\in \O}|\lt^-
 (x,t)|\intr (\lt^2+\lt\ltt +\ltt^2) dx,
 \eq
  where we used (\ref{basic}) again.
  Similarly to the above, combining this with
(\ref{2.1}) and (\ref{2.8a}), yields
 \bb
 \frac{d}{dt}\intr |\o (x,t)|^2 dx\geq \left[\inf_{x\in \O}\lt^+ (x,t)
-2\sup_{x\in \O}|\lt^-
 (x,t)|\right]\intr |\o (x,t)|^2 dx,
 \ee
 and, applying the Gronwall lemma we finish the first inequality of
 (\ref{important}).
 $\square$\\

\begin{cor}
Let $v_0 \in \mathbb{H}^m_\s $ be given, and  $\lo (x,t),\lt
(x,t),\ltt (x,t)$ are as in Theorem 2.2. Suppose
 \bb\label{check}
\lim\sup_{t\to T_*} \|\o (t)\|_{L^2} =\infty .
  \ee
Then, necessarily
 \bb \int_0 ^{T_*} \|\lt ^+
(t)\|_{L^\infty} dt =\infty.
 \ee
\end{cor}
{\bf Proof.} We just observe that from (\ref{important}), we have
immediately
$$
 \|\o(t)\|_{L^2}\leq \|\o_0\|_{L^2} \exp \left(\int_0 ^t \|\lt^+
 (s)\|_{L^\infty}ds\right).
 $$
 This implies the corollary.
$\square$\\

\noindent{\bf Remark 2.3} The above corollary says that if
singularity happens in the $L^2$ norm of vorticity, then it should
be caused by the uncontrollable
intensification of planar stretching.\\

\noindent{\bf Remark 2.4} In the 3D Navier-Stokes equations the
$L^2$ norm singularity of vorticity is equivalent to the that of any
high norms in Sobolev space(see e.g. \cite{tem1, con3}). Hence, the
above corollary says that the regularity/singularity of the 3D
Navier-Stokes equations are controlled by the integral, $
 \int_0 ^{t} \|\lt ^+
(s)\|_{L^\infty} ds . $\\

\noindent{\bf Remark 2.5} In \cite{he} the author also investigated
another sufficient condition for the singularity of $L^2$ norm of
vorticity of the 3D Euler equations, using simultaneously the
eigenvector and eigenvalues of the deformation tensor and the
hessian of the pressure. Our condition is completely different from
it, and has direct physical interpretation.\\

\section{Applications for some classes of initial data}
\setcounter{equation}{0}

 In order to state our theorem in this section we introduce
some definitions. Given a differentiable vector field $f=(f_1 ,f_2
,f_3 )$ on $\Bbb T^3$, we denote by the scalar field $\l_i (f)$,
i=1,2,3, the eigenvalues of the deformation tensor associated with
$f$. Below we always assume the ordering, $ \lo (f)\geq \lt (f)\geq
\ltt (f). $ We also fix $m>5/2$ below. We recall that if $f\in
\mathbb{H}^m _\s $, then $\lo (f)+\lt (f)+\ltt (f)=0$, which is
another representation of div $f=0$.

Let us begin with introduction of  admissible classes
$\mathcal{A}_\pm$ defined by
$$\mathcal{A}_+=\{ f\in \mathbb{H}^m _\s (\Bbb T^3)\, | \, \inf_{x\in \O}\lt (f)(x)
>0 \,\},$$
and
$$\mathcal{A}_-=\{ f\in \mathbb{H}^m_\s (\Bbb T^3)\, |
\sup_{x\in \O}\lt (f)(x)<0 \, \}.$$
 Physically $\mathcal{A}_+$ consists of solenoidal vector fields
 with planar stretching(see Remark 2.2) everywhere, while $\mathcal{A}_-$ consists
 of everywhere linear stretching vector fields. Although they do not represent
 real physical flows, they might be useful in the study of searching
 initial data leading to finite time singularity for the 3D Euler
 equations.

 Given $v_0
\in \mathbb{H}_\s ^m$, let $T_*(v_0)$ be the maximal time of unique
existence of solution in $\mathbb{H}_\s ^m$ for the system
(1.1)-(1.3).
 Let $S_t : \mathbb{H}^m_\s \to
\mathbb{H}^m_\s$ be the solution operator, mapping from initial data
to the
 solution $v(t)$.
Given $f\in \mathcal{A}_+$, we define the first zero touching time
of $\lt (f)$ as
$$ T(f)=
\inf\{ t\in (0, T_* (v_0)) \, | \, \mbox{$\exists x\in \Bbb T^3 $
such that $\lt (S_t f ) (x)<0 $}\}.
$$
Similarly for $f\in \mathcal{A}_-$, we define
$$ T(f)=
\inf\{ t\in (0, T_* (v_0)) \, | \, \mbox{$\exists x\in \Bbb T^3 $
such that $\lt (S_t f ) (x)>0 $}\}.
$$
The following theorem is actually an immediate corollary of Theorem
2.2, combined with the above definition of $\mathcal{A}_\pm$ and
$T(f)$. We just observe that for $v_0 \in \mathcal{A}_+ $(resp.
$\mathcal{A}_- $) we have $\lt^-=0, \lt^+=\lt$(resp. $\lt^+ =0,
\lt^-=\lt$) on $\O \times (0,T( v_0 ))$.

\begin{theorem}
Let $v_0\in \mathcal{A}_\pm$ be given. We set $\lo(x,t)\geq
\lt(x,t)\geq \ltt(x,t)$ as the eigenvalues of the deformation tensor
associated with $v(x,t)=(S_t v_0)(x)$ defined $t\in (0, T(v_0 ))$.
Then, for all $t\in (0, T(v_0 ))$ we have the following
estimates:\\
 (i) If $v_0 \in \mathcal{A}_+$, then
 \bb\label{a}
 \exp\left( \frac12\int_0 ^t
 \inf_{x\in\O} |\lt (x,s)| ds\right)\leq \frac{\|\o (t)\|_{L^2}}{ \|\o_0\|_{L^2}}
  \leq  \exp\left( \int_0 ^t
 \sup_{x\in\O} |\lt (x,s)| ds\right) .
 \ee
(ii) If $v_0 \in \mathcal{A}_-$, then
 \bb\label{b}
 \exp\left( - \int_0 ^t
 \sup_{x\in\O} |\lt (x,s)| ds\right)\leq \frac{\|\o (t)\|_{L^2}}{ \|\o_0\|_{L^2}}
  \leq  \exp\left( -\frac12 \int_0 ^t
 \inf_{x\in\O} |\lt (x,s)| ds\right).
 \ee
\end{theorem}
{\bf Remark 3.1} If we have the comparability conditions,
 \bqn
\inf_{x\in\O} |\lt (x,t)|&\simeq& \sup_{x\in\O} |\lt (x,t)|\simeq
g(t) \quad \forall t\in (0, T(v_0 )), \eqn
 which is the case for sufficiently small box $\Bbb T^3$, then we have
 \bqn
 \frac{\|\o (t)\|_{L^2}}{ \|\o_0\|_{L^2}} \simeq\left\{
 \aligned &\exp\left(\int_0 ^t g(s)ds \right) \quad \mbox{if}\quad v_0 \in
 \mathcal{A}_+\\
 &\exp\left(-\int_0 ^t g(s)ds \right) \quad \mbox{if}\quad v_0
\in
 \mathcal{A}_- \endaligned \right.
 \eqn
 for $t\in (0, T(v_0 ))$.
In particular, if we could find $v_0 \in \mathcal{A}_+$ such that
\bb \label{sing} \inf_{x\in\O} |\lt (x,t)| \gtrsim
O\left(\frac{1}{t_*-t}\right)
 \ee for time
interval near $t_*$, then such data would lead to singularity at
$t_*$.

Below we have some decay in time estimates for some ratio of
eigenvalues.

\begin{theorem}
Let $v_0 \in \mathcal{A}_\pm$ be given, and we set $\lo(x,t)\geq
\lt(x,t)\geq \ltt(x,t)$ as  in Theorem 3.1.
  We define
 \bb\label{3.1}
 \e (x,t)=\frac{|\lt (x,t) |}{\l (x,t)}\quad \forall (x,t)\in
\Bbb T^3 \times (0, T(v_0)),
  \ee where we set
$$\l (x,t)=\left\{\aligned \l _1 (x,t)
\quad \mbox{ if} \quad v_0 \in  \mathcal{A}_+\\
 -\l _3 (x,t) \quad \mbox{ if} \quad v_0 \in  \mathcal{A}_- .
\endaligned
\right.
$$
  Then, there exists a constant $C=C(v_0, |\Om|)$, with $|\Om|$ denoting the volume
  of the box $\Om=\Bbb T^3$, such that
 \bb \label{3.2}
 \inf_{(x,s)\in \Bbb T^3\times (0,t)} \e (x,s)< \frac{C}{\sqrt{t}}\quad
 \forall t\in (0, T(v_0)).
 \ee
 \end{theorem}
{\bf Remark 3.2} Regarding the problem of searching  a finite time
 blowing up solution, again,  the proof of
 the above theorem, in particular, the
estimate (\ref{3.7}) below,
 combined with Remark 2.3, suggests the following,
  :\\
 Given $\delta >0$, let us suppose we could find $v_0 \in
 \mathcal{A}_+$ such that for the associated solution $v(x,t)=(S_t v_0)(x)$
 the estimate
 \bb\label{3.3}
 \inf_{(x,s)\in \Bbb T^3 \times (0,t)}\e (x,s)\gtrsim
 O\left(\frac{1}{t^{\frac12+\delta}}\right),
 \ee
 holds true, for sufficiently
 large time $t$.
Then such $v_0$ will lead to the finite time singularity. In order
to check the  behavior (\ref{3.3}) for a given solution we need a
sharper and/or localized
 version of the equation (\ref{2.1}) for the dynamics of eigenvalues
 of the deformation tensor.\\
\ \\
\noindent{\bf Proof of  Theorem 3.2} We divide the proof  the into
two separate cases.\\

 \noindent \underline{{(i) The case  $v_0 \in \mathcal{A}_+$:}}\\

 We parameterize the eigenvalues of the deformation
 tensor of the solution $v(x,t)$
 of (1.1)-(1.3) by
 $$ \lo (x,t)= \l (x,t) >0 , \,\lt=\e (x,t)\l (x,t) >0, \, \ltt
 (x,t)=-(1+\e (x,t))\l (x,t) <0.$$
 for all $(x,t)\in \Bbb T^3 \times (0,T)$. We observe that
 \bb
 \label{3.4}
 0<\e (x,t)\leq 1 \quad \forall (x,t)\in \Bbb T^3 \times [0, T(v_0)).
 \ee
 The equation (\ref{2.1}) can be written as
 \bb\label{3.5}
 \frac{d}{dt} \intr \l ^2 (\e^2 +\e +1) dx =2 \intr \l^3 (\e^2 +\e)
 dx \quad \forall t\in (0, T(v_0)).
 \ee
From the estimate
 \bqn
 \lefteqn{\intr  \l ^2 (\e^2 +\e +1) dx =\intr \l^2 (\e^2 +\e)^{\frac23}
 \frac{(\e^2 +\e +1)}{ (\e^2 +\e)^{\frac23}} dx}\hspace{.2in}\\
 &\leq&\left[ \intr \l^3 (\e^2 +\e)
 dx\right]^{\frac23} \left[\intr \frac{(\e^2 +\e +1)^3}{ (\e^2
 +\e)^{2}}dx\right]^{\frac13}\\
&\leq&\frac{3}{\sqrt[3]{4}}\left[ \intr \frac{1}{\e
^4}dx\right]^{\frac13}  \left[ \intr \l^3 (\e^2 +\e)
 dx\right]^{\frac23},
\eqn
 where we used (\ref{3.4}),
 we have
 $$
 \intr \l^3 (\e^2 +\e)dx \geq
 \frac{2}{\sqrt{27}}\left[ \intr \frac{1}{\e
^4}dx\right]^{-\frac12}\left[\intr \l ^2 (\e^2 +\e +1)
dx\right]^{\frac32}
 dx,
 $$
 which, combined with (\ref{3.5}), yields
 \bb\label{3.6}
 \frac{d}{dt} \intr
 \l ^2 (\e^2 +\e +1) dx \geq  \frac{4}{\sqrt{27}}\left[ \intr \frac{1}{\e
^4}dx\right]^{-\frac12} \left[\intr \l ^2 (\e^2 +\e +1)
dx\right]^{\frac32}.
 \ee
  Setting
  $$ y(t)=\left[\intr \l ^2 (\e^2 +\e +1) dx\right]^{\frac12}, $$
  we have
  $$\frac{dy}{dt}\geq  \frac{2}{\sqrt{27}}\left[ \intr \frac{1}{\e
^4}dx\right]^{-\frac12}y^2.$$
  Solving the differential inequality, we have
  $$
  y(t)\geq \frac{y_0}{1-\frac{2y_0}{\sqrt{27}}\int_0 ^t\left[ \intr \frac{1}{\e
^4}dx\right]^{-\frac12} ds }.
  $$
  Since
  $y^2(t)= \frac{1}{2} \|\o (t)\|_{L^2}^2$,
  we have just derived
  $$
  \|\o (t)\|_{L^2} \geq \frac{\sqrt{2}
  \|\o _0\|_{L^2}}{\sqrt{2}-\frac{2\|\o_0\|_{L^2}}{\sqrt{27}}
  \int_0 ^t\left[ \intr \frac{1}{\e
^4}dx\right]^{-\frac12} ds }\quad  \forall t\in [0, T(v_0)).
  $$
  Since the denominator should be positive for all $t\in [0,
  T(v_0)]$, we obtain that
  $$
\frac{2\|\o_0\|_{L^2}}{\sqrt{27}}
  \int_0 ^t\left[ \intr \frac{1}{\e
^4}dx\right]^{-\frac12} ds  <\sqrt{2}.
 $$
 Estimating from below the left hand side, we are lead to the
 inequality
  \bb\label{3.7}
 t \inf_{(x,s)\in \Bbb T^3 \times (0,t)} \e ^2(x,s)
 \leq |\Omega |^{\frac12}\int_0 ^t\left[ \intr \frac{1}{\e
^4}dx\right]^{-\frac12} ds  < \frac{\sqrt{27}|\Omega |^{\frac12}
}{\sqrt{2}\|\o_0\|_{L^2}},
 \ee
 which implies (\ref{3.2}) for the case $v_0 \in \mathcal{A}_+$.\\
 \ \\
\noindent \underline{{(ii) The case $v_0 \in \mathcal{A}_-$:}}\\

 In this case parameterize the eigenvalues as
 $$
 \lo (x,t)= (1+\e (x,t))\l (x,t) >0 , \quad\lt=-\e (x,t)\l (x,t) >0
 $$
 $$
\ltt
 (x,t)=-\l (x,t)>0,
$$
where as previously we have $0<\e (x,t) \leq 1 $ for all $(x,t)\in
\O \times (0, T(v_0 ))$.
 The equation (\ref{2.1}) can now be written as
 \bb\label{3.8}
 \frac{d}{dt} \intr \l ^2 (\e^2 +\e +1) dx =-2 \intr \l^3 (\e^2 +\e)
 dx .
 \ee
Similarly to the above, we  obtain
 \bb\label{3.9}
 \frac{d}{dt} \intr
 \l ^2 (\e^2 +\e +1) dx \leq -\frac{2}{\sqrt{27}}\left[ \intr \frac{1}{\e
^4}dx\right]^{-\frac12}\left[\intr \l ^2 (\e^2 +\e +1)
dx\right]^{\frac32}.
 \ee
  Hence, by similar procedure to the previous case, we have
  \bb\label{3.10}
  \|\o (t)\|_{L^2} \leq \frac{\sqrt{2}
  \|\o _0\|_{L^2}}{\sqrt{2}+\frac{2\|\o_0\|_{L^2}}{\sqrt{27}}
  \int_0 ^t\left[ \intr \frac{1}{\e
^4}dx\right]^{-\frac12} ds }.
  \ee
Now we recall the helicity conservation(see e.g. \cite{maj}),
$$
 H(t) =\intt v(x,t)\cdot \o(x,t) dx =\intt v_0 (x)\cdot \o _0 (x)dx=H_0,
 $$
 which implies
 \bb\label{3.11}
 H_0 \leq \|v(t)\|_{L^2} \|\o (t)\|_{L^2}=\|v_0\|_{L^2} \|\o
 (t)\|_{L^2}= \sqrt{2 E_0}\|\o
 (t)\|_{L^2} ,
\ee
 where we used the energy conservation
 $$
 E(t)=\frac12 \|v(t)\|_{L^2}^2 =\frac12 \|v_0\|_{L^2}^2=E_0.
 $$
 Combining (\ref{3.10}) with (\ref{3.11}), we have
$$
\frac{H_0}{\sqrt{2E_0}} \leq \frac{\sqrt{2}
  \|\o _0\|_{L^2}}{\sqrt{2}+\frac{2\|\o_0\|_{L^2}}{\sqrt{27}}
  \int_0 ^t\left[ \intr \frac{1}{\e
^4}dx\right]^{-\frac12} ds },
$$
from which we  derive
 \bb\label{3.12} \int_0 ^t \left[ \intt  \frac{1}{\e
^4}dx
 \right]^{-\frac12} ds \leq\sqrt{27} \left(\frac{\sqrt
 E_0}{H_0} -\frac{1}{\sqrt{2}\|\o_0\|_{L^2}}\right).
\ee
 Estimating from below the left hand side of (\ref{3.12}), we deduce
 \bb
 t\inf_{(x,s)\in \Bbb T^3 \times [0,s] } \e^2 (x,s)\leq
\sqrt{27} |\Om|^{\frac12}\left(\frac{\sqrt
 E_0}{H_0} -\frac{1}{\sqrt{2}\|\o_0\|_{L^2}}\right)
 \ee
 for all $t\in (0, T(v_0 ))$.
 This finishes the proof of (\ref{3.2}) for $v_0 \in \mathcal{A}_-$.
$\square$ \\
\ \\
$$ \mbox{\bf Acknowledgements} $$
 This work was supported by Korea Research Foundation
Grant KRF-2002-015-CS0003.

\end{document}